\magnification=\magstep1   
\input amstex
\UseAMSsymbols
\input pictex
\vsize=23truecm
\NoBlackBoxes
\parindent=18pt
  
   \font\rmk=cmr8      \font\ttk=cmtt8

\font\gross=cmbx10 scaled\magstep1 

\def\op{{\text{op}}}
\def\mod{\operatorname{mod}}

\def\Hom{\operatorname{Hom}}
\def\End{\operatorname{End}}
\def\Ann{\operatorname{Ann}}

\def\rad{\operatorname{rad}}
\def\add{\operatorname{add}}

\def\soc{\operatorname{soc}}
\def\top{\operatorname{top}}

\def\bdim{\operatorname{\bold{dim}}}

\def\Rahmen#1%
   {$$\vbox{\hrule\hbox%
                  {\vrule%
                       \hskip0.5cm%
                            \vbox{\vskip0.3cm\relax%
                               \hbox{$\displaystyle{#1}$}%
                                  \vskip0.3cm}%
                       \hskip0.5cm%
                  \vrule}%
           \hrule}$$}

\vglue1.8truecm

\centerline{\gross Quiver Grassmannians and 
Auslander varieties}
	\medskip
\centerline{\gross 
for wild algebras.}
	\bigskip
\centerline{Claus Michael Ringel}
	\bigskip\bigskip
{\bf Abstract.} Let $k$ be an algebraically closed field and 
$\Lambda$ a finite-dimensional $k$-algebra. Given a $\Lambda$-module $M$, 
the set $\Bbb G_{\bold e}(M)$ of all submodules of
$M$ with dimension vector $\bold e$ is called a quiver Grassmannian. 
If $D,Y$ are $\Lambda$-modules, then we
consider $\Hom(D,Y)$ as a $\Gamma(D)$-module, where $\Gamma(D) = 
\End(D)^\op$, and the Auslander varieties for $\Lambda$ are the
quiver Grassmannians of the form $\Bbb G_{\bold e}\Hom(D,Y)$. 
Quiver Grassmannians, 
thus also Auslander varieties 
are projective varieties and it is known that every projective 
variety occurs in this way. There is a tendency to relate this fact
to the wildness of quiver representations and the aim of this note is to clarify
these thoughts: 
We show that for an algebra $\Lambda$ which is (controlled) wild, 
any projective
variety can be realized as an Auslander variety, but not necessarily
as a quiver Grassmannian. 

	\bigskip\bigskip
{\bf 1. Introduction.} 
	\medskip 
Let $k$ be an algebraically closed field and 
$\Lambda$ a finite-dimensional 
$k$-algebra.  A {\it dimension vector} $\bold d$ for $\Lambda$
is a function defined on the set of isomorphism classes of simple $\Lambda$-modules $S$ with 
values $d_S$ being non-negative integers. If $M$ is a $\Lambda$-module, its dimension vector $\bdim M$
attaches to the simple module $S$ the Jordan-H\"older multiplicity $(\bdim M)_S = [M:S].$ 

Given a $\Lambda$-module $M$, 
the set $\Bbb G_{\bold e}(M)$ of all submodules of
$M$ with dimension vector $\bold e$ is called a quiver Grassmannian. 
Quiver Grassmannians are projective varieties and every projective 
variety occurs in this way (see the Appendix). 
If $D,Y$ are $\Lambda$-modules, then we
consider $\Hom(D,Y)$ as a $\Gamma(D)$-module, where $\Gamma(D) = 
\End(D)^\op$. The easiest way to define the Auslander varieties 
for $\Lambda$ is to say that they are just the
quiver Grassmannians $\Bbb G_{\bold e}\Hom(D,Y)$ (here, we rely on the
Auslander bijections; 
the proper definition of the Auslander varieties
would have to refer to right equivalence classes of right $D$-determined maps 
ending in $Y$, see [Ri]). The Auslander varieties are part of
Auslander's approach of describing the global directedness of the category
$\mod\Lambda$. Let as add that the quiver Grassmannians for $\Lambda$ are special 
Auslander varieties, namely the 
Auslander varieties $\Bbb G_{\bold e}\Hom(D,Y)$ with $D = \Lambda$.

According to Drozd [D1], any finite dimensional $k$-algebra is either tame
or wild 
(note that 
there are few tame algebras, most of the algebras are wild; for example,
the path algebra of a connected quiver is tame only in case we deal with
a Dynkin or an extended Dynkin quiver). It has been  conjectured that wild
algebras are actually controlled wild (the definition will be recalled in section 2).
A proof of this conjecture has been
announced by Drozd [D2] in 2007, but apparently it has not yet been published. 
We show that for a fixed (controlled) wild algebra $\Lambda$, any projective
variety can be realized as an Auslander variety, but not necessarily
as a quiver Grassmannian. 
	\bigskip
{\bf 2. Controlled wild algebras.}
	\medskip 
We denote by $\mod\Lambda$ the category of all (finite-dimensional left)
$\Lambda$-modules. Let $\rad$ be the radical of $\mod\Lambda$, this is the ideal
generated by all non-invertible maps between indecomposable modules.
If $\Cal C$ is a collection of object
of $\mod \Lambda$, we denote by $\add \Cal C$ the closure under direct sums and direct summands.
For every pair $X,Y$ of $\Lambda$-modules, 
$\Hom(X,\Cal C,Y)$ denotes the subgroup of $\Hom(X,Y)$ 
given by the maps $X \to Y$ which factor through
a module in $\add\Cal C$. 
	\medskip
Here is now the definition.
The algebra $\Lambda$ is said to be {\it controlled wild} provided for any
finite-dimensional $k$-algebra $\Gamma$, there is an exact embedding
functor $F\:\mod \Gamma \to \mod \Lambda$ and a full subcategory $\Cal C$ of
$\mod \Lambda$ (called the {\it control class}) 
such that for all $\Gamma$-modules $X,Y$,
the subgroup $\Hom(FX,\Cal C,FY)$ is contained in $\rad (FX,FY)$ and we have
$$
 \Hom(FX,FY) = F\Hom(X,Y) \oplus \Hom(FX,\Cal C,FY).
$$

In order to check that $\Lambda$ is controlled wild, it is sufficient to
exhibit such a functor $F$ for just one suitable algebra $\Gamma$, for example for
the 3-Kronecker algebra (this is the path algebra of the quiver
with two vertices, say $a$ and $b,$ and three arrows $b \to a$). 

We also mention that $\Lambda$ is said to be {\it strictly wild} 
provided for any
finite-dimensional $k$-algebra $\Gamma$, there is a full exact embedding
functor $F\:\mod \Gamma \to \mod \Lambda$ (thus, strictly wild algebras are
controlled wild and we can take as control class $\Cal C$ the zero subcategory).
The 3-Kronecker algebra is a typical strictly wild algebra. 
	\bigskip
{\bf 3. Auslander varieties}
		\medskip 
{\bf Proposition 1.} {\it Let $\Lambda$ be a finite-dimensional 
$k$-algebra which is controlled wild. Let $V$ be any projective variety. 
Then there
are $\Lambda$-modules $D,Y$ and a dimension vector $\bold e$ for $\Gamma(D)$ 
such that $\Bbb G_{\bold e} \Hom(D,Y)$ is of the form $V$.}
	\medskip
The special case of a strictly wild algebras has been considered already in [Ri]. 
The proof of Proposition 1 will be given in this section. We start with the following Lemma.
	\medskip 
{\bf Lemma 1.} {\it Given finitely many modules $X_1,\dots, X_n$ and a set $\Cal C$ of modules in $\mod\Lambda$,
then there is a module $C$ in $\add\Cal C$ such that $\Hom(X_i,\Cal C,X_j) = \Hom(X_i,C,X_j)$ for all $i,j.$}
	\medskip
Proof: Let $X = \bigoplus X_i$. It is sufficient to show that $\Hom(X,C,X) = \Hom(X,\Cal C,X)$ for
some module $C\in \add \Cal C$. 
Since  the subgroups
$\Hom(X,C,X)$ with $C \in \add\Cal C$ are subspaces
of the  finite-dimensional vector space $\Hom(X,X)$, 
there is  $C\in \add\Cal C$
such that $\Hom(X,C,X)$ is of maximal dimension. 
Let $C'\in \add\Cal C$. Then also $C\oplus C'$ belongs to
$\add\Cal C$ and we have $\Hom(X,C,X) \subseteq
\Hom(X,C\oplus C',X)$. The maximality of the dimension of $\Hom(X,C,X)$ implies 
that $\Hom(X,C,X) =
\Hom(X,C\oplus C',X),$ and thus $\Hom(X,C',X) \subseteq \Hom(X,C,X)$.
But $\Hom(X,\Cal C,X) = \bigcup_{C'} \Hom(X,C',X)$.
	
	\bigskip 
{\bf Lemma 2.} {\it Let $R$ be a finite-dimensional algebra and $e$ an idempotent
in $R$. Let $N$ be an $R$-module and $\bold c = \bdim ReN.$ If $\bold g$ is a 
dimension vector and $U$ belongs to
$\Bbb G_{\bold g+\bold c} N$, then $U \supseteq ReN$, and $U/ReN$ is an element of
$\Bbb G_{\bold g}(N/ReN).$ As a consequence, the varieties $\Bbb G_{\bold g+\bold c} N$ 
and $\Bbb G_{\bold g}(N/ReN)$
can be identified.}

 	\medskip 
Proof. Let $U$ be an element of $\Bbb G_{\bold g+\bold c} N$. We want to show that 
$U \supseteq ReN$. Given dimension vectors $\bold d, \bold d'$ for $\Lambda$, one writes
$\bold d'\le \bold d$ provided $\bold d-\bold d'$ has non-negative coefficients.
Since $\bdim U = \bold g+\bold c,$ we have $\bdim U \ge \bold c.$ Let $S$ be a simple
$R$-module with $eS \neq 0$. Then
$$
 [U:S] = (\bdim U)_S \ge \bold c_S = (\bdim \Lambda eN)_S = [\Lambda eN:S],
$$
and therefore $eN \subseteq U$, thus also $\Lambda eN \subseteq U.$ 
	\bigskip 

{\bf Proof of proposition 1.} Let $V$ be a projective variety. There is a finite-dimensional
algebra $\Gamma$, a $\Gamma$-module $M$ and a dimension vector $\bold g$ for $\Gamma$
such that $\Bbb G_{\bold g}M$ is of the form $V$ (see the Appendix). Since
$\Lambda$ is controlled wild, there is a controlled embedding $F$ of $\mod
\Gamma$ into $\mod\Lambda$, say with control class $\Cal C$.
Let $G = F({}_\Gamma\Gamma)$ and $Y = F(M).$ According to Lemma 1, there is
$C\in \add\Cal C$ such that $\Hom(G,C,G) = \Hom(G,\Cal C,G)$ and 
$\Hom(G,C,Y) = \Hom(G,\Cal C,Y)$. Let $D = G\oplus C$
and $R = \End(D)^\op.$ Let $e_G$ be the projection of $D$ onto $G$ with kernel $C$
and $e = e_C$ the projection of $D$ onto $C$ with
kernel $G$, both $e_G, e_C$ considered as elements of $R$. Note that
$$
\align
 \Hom(D,D) &= \Hom(G\oplus C,G\oplus C) = \Hom(G,G)\oplus \Hom(G,C) \oplus \Hom(C,G) \oplus \Hom(C,C)\cr
   &= F(\Hom(\Gamma,\Gamma)) 
 \oplus \Hom(G\oplus C,\Cal C,G\oplus C) \oplus \Hom(G,C) \oplus \Hom(C,G) \oplus \Hom(C,C),
\endalign
$$
and 
$$
 e\Hom(D,D)e = \Hom(G\oplus C,\Cal C,G\oplus C) \oplus \Hom(G,C) \oplus \Hom(C,G) \oplus \Hom(C,C).
$$
It follows that the map $\gamma \mapsto F(\gamma)\in 
e_GRe_G$ yields an isomorphism $\Gamma \to R/ReR.$ 
Also, we are interested in the $R$-module $N = \Hom(G\oplus C,Y).$ Here, we have 
$$
\align
 N = \Hom(G\oplus C,Y) &= \Hom(G\oplus 0,Y) \oplus \Hom(0\oplus C,Y) \cr
                   &= F\Hom(\Gamma,M) \oplus \Hom(G\oplus 0,C,Y) \oplus 
                      \Hom(0\oplus C,Y).
\endalign
$$
If we multiply $N$ with the element $e = e_C\in R$, we obtain
$$
 eN =  \Hom(0\oplus C,Y),
$$
thus 
$$
 ReN = 
  R\Hom(0\oplus C,Y) = \Hom(G\oplus 0,C,Y) \oplus 
                      \Hom(0\oplus C,Y).
$$
This shows that $N/ReN$ is canonically isomorphic to $F\Hom(\Gamma,M)$ as an
$R$-module. Of course, these modules are annihilated by $e$, thus they are
$R/ReR$-modules and as we know $R/ReR = \Gamma.$ 

It remains to apply Lemma 2. 
	\bigskip\bigskip
\vfill\eject
{\bf 4. Quiver Grassmannians.}
	\medskip 
{\bf Proposition 2.} {\it There are controlled wild algebras $\Lambda$ such 
that not every projective variety can be realized as a
quiver Grassmannian of a $\Lambda$-module.}
	\medskip 
Proof. Let $\Lambda$ be any local radical square zero $k$-algebra of dimension at least 4
(thus $\Lambda = k[T_1,\dots,T_n]/(T_1,\dots,T_n)^2$ with $n \ge 3$). 
It is well-known (and easy to see) that such an algebra 
is controlled wild. 
Let $M$
be a $\Lambda$-module. The Grothendieck group $K_0(\Lambda)$ is free of rank one,
thus the quiver Grassmannians are of the form $\Bbb G_i(M)$ with $i$ a non-negative
integer (the elements of $\Bbb G_i(M)$ are the submodules of $M$ of dimension $i$).
In order to determine the possible varieties $\Bbb G_i(M)$, we can assume that $i 
\le\dim \soc M$. Namely, if $i > \dim\soc M,$ then we consider the dual module
$M^*$ and the quiver Grassmannian $\Bbb G_{d-i}(M^*)$, where $d = \dim M = \dim M^*$.
On the one hand, the varieties $G_i(M)$ and $G_{d-i}(M^*)$ are obviously isomorphic,
on the other hand we have $d-i < \dim M - \dim\soc M \le \dim M - \dim\rad M = 
\dim \top M = \dim\soc(M^*)$, here we have used that $\rad M \subseteq \soc M$.  

Thus, let $i \le s = \dim\soc M.$ The submodules $U$ of $\soc M$ of dimension $i$ 
are just the subspaces of $\soc M$ of dimension $i$
considered as a $k$-space, thus they form
the usual Grassmannian $\Bbb G_i(\soc M) = \Bbb G_i(k^s)$, in particular, this
is an irreducible (and rational) variety. 
Now let $U$ be a submodule of $M$ of dimension $i$
which is not contained in $\soc M$. We claim that there is a projective line 
in $\Bbb G_i(M)$ which contains both $U$ and a submodule $U'$ of $\soc M$.
Namely, let $b_1,\dots, b_t$ be a basis
of $\soc U$, and extend it to a basis $b_1,\dots,b_i$ of $U$. Now $b_1,\dots, b_t$
are linearly independent elements of the socle of $M$. Since $i \le \dim\soc M$,
there is an $i$-dimensional subspace $U'$ inside $\soc M$ which contains
$\soc U' = U \cap \soc M.$ 
We can extend the basis $b_1,\dots, b_t$ 
of $\soc U$ to a basis of $U'$, say $b_1,\dots, b_t, b'_{t+1},\dots, b'_i$,
For $\lambda = (\lambda_0:\lambda_1) \in \Bbb P^1$, define $U_\lambda$
as the subspace of $M$ with basis the elements $b_1,\dots, b_t$ as well as the
elements
$\lambda_0b_j+\lambda_1b'_j$ where $t < j \le i.$ Of course, $U_\lambda$ is a
submodule of $M$. On the one hand, we have $U_{(1:0)} = U$, on the other hand, 
$U_{(0:1)} = U'$ is an $i$-dimensional
submodule of $M$ which lies inside the socle of $M$. 

Since the Grassmannian $\Bbb G_i(\soc M)$ is (rational and) connected and for any element $U\in 
\Bbb G_i(M)$ there is a $\Bbb P^1$-family of submodules which contains $U$ and
an element in $\Bbb G_i(\soc M)$, it follows that also $\Bbb G_i(M)$
is connected (even rationally connected, see [Ha]).

	\bigskip\bigskip
{\bf 5. Open questions.}
	\medskip
We have shown that any projective variety occurs as an Auslander variety for
any (controlled) wild algebra. It seems that the Auslander varieties for the tame
algebras are quite restrictive --- is there a special property which all have?
Such a result would provide a characterization of the tame-wild dichotomy in
terms of Auslander varieties. 

We have shown that dealing with local algebras with radical square zero, all
quiver Grassmannians are rationally connected. Are there further properties which they share?
On the other hand, the quiver Grassmannians for wild algebras should be of quite
a general nature. Is there a class of varieties which can be realized as 
quiver Grassmannians for all wild algebras, but not for tame algebras?
 
	\bigskip\bigskip 
\vfill\eject
{\bf Appendix. Every projective variety is a quiver Grassmannian.}
	\medskip
We say that a $k$-module $M$ is a {\it brick} provided $\End(M) = k.$
	\medskip 
{\bf Proposition (Reineke).} {\it 
Every projective variety is a quiver Grassmannian $\Bbb G_{\bold e}M$ for
a brick $M$.}
	\medskip
Let us outline a proof, following Van den Bergh [L]. Let $V$ be a projective
variety. We can assume $V$ is a closed subset of the projective space $\Bbb P^n$,
defined by the vanishing of homogeneous polynomials $f_1,\dots,f_m$ of degree 2.
Let $\Delta$ be the Beilinson quiver with 3 vertices, say $a,b,c$, with $n+1$ arrows
$b\to a$ labeled $x_0,\dots,x_n$ as well as $n+1$ arrows $c\to b$, also
labeled $x_0,\dots,x_n$. The path algebra of $\Delta$ with
all the relations $x_ix_j = x_jx_i$ (whenever this
makes sense) is called the Beilinson algebra. Let $\Lambda$ be the factor
algebra of this Beilinson algebra taking the elements $f_1,\dots,f_m$
as additional relations (obviously, these elements may be considered as 
linear combinations of paths of length 2 in the quiver $\Delta$). 
Take $M = I_{\Lambda}(a),$ the indecomposable injective $\Lambda$-module 
corresponding to the vertex $a$, and take $\bold e = (1,1,1)$. Note that
the elements of $\Bbb G_{\bold e}M$ are just all the serial $\Lambda$-modules,
one from each isomorphism class (we call a module $X$ {\it serial}, provided
it has a unique composition series).
	\medskip
Here are some remarks on the history: 
The title of the appendix is also the title of a recent paper [Re] by Reineke,
who answered in this way a question by Keller. 
The 2-page paper attracted a lot of interest, see for example 
blogs by Le Bruyn [L] (with the proof by Van den Bergh presented above) and
by Baez [Bz]. Actually, the construction given in the proof of Van den Bergh
is much older, it has been used before by several 
mathematicians dealing with related problems and may, of course, be traced
back to Beilinson [Be].
	 
The quiver Grassmannians
play an important role in the representation theoretical approach to
cluster algebras. Here one deals with the quiver Grassmannians $\Bbb G_{\bold e}M$,
where $M$ is a quiver representation without self-extensions. It has been
asserted by Caldero and Reineke [CR] that the quiver Grassmannians $\Bbb G_{\bold e}M$,
where $M$ is a quiver representation without self-extensions, are very special:
{\it if $\Bbb G_{\bold e}M$ is non-empty, then the Euler characteristic of $\Bbb G_{\bold e}M$
is positive}; a complete proof was given by Nakajima [N] and Qin [Q]. 
Note that an indecomposable quiver representation without 
self-extensions is a brick. If we consider the bricks $M$ constructed by  
Van den Bergh as representations of the quiver $\Delta$, then such a $k\Delta$-module 
will have self-extensions, but it seems to be remarkable to observe that $M$ has 
no self-extensions when it is considered as a $k\Delta/\Ann(M)$-module, 
where $\Ann(M)$ is the annihilator of $M$ in $k\Delta$
(in contrast, the examples constructed by Reineke are usually faithful quiver 
representations). Namely, we have $k\Delta/\Ann(M) = \Lambda$, and by construction, 
$M$ is an injective $\Lambda$-module.
We should recall that for any ring $R$, an $R$-module is said to be {\it quasi-injective}
provided for any submodule $U$ of $M$ any map $U \to M$ can be extended
to an endomorphism of $M$; for an artinian ring, a module $M$ is quasi-injective if and
only if $M$ considered as an $R/\Ann(M)$-module is injective. 
	
There is a tendency to relate the fact that every variety is a quiver Grassmannian to the
tame-wild dichotomy as established by Drozd [D1]. For example, Baez [B] writes that one may suppose that
this is {\it just another indication of the 'wildness' of quiver representations 
once we leave the safe waters of Gabriel's theorem.} The aim of this note was to clarify
these thoughts.
	 
Some mathematicians (see [L],[V]) refer in this context to  ``Murphy's law'':
{\it Anything that can go wrong, will go wrong}
(as formulated in the Wikipedia [W]), or: {\it Anything that can happen, will happen.} 
But one should be aware that Murphy's law
may be a challenging assertion in daily life, 
but it is just a tautology when we consider mathematical questions. Indeed,
in mathematics, if we know (that means: if we can prove) that something 
{\bf does not} happen, then of course we have a proof that it {\bf cannot} happen. 
	\bigskip\bigskip 
{\bf References.}
	\medskip
\item{[Bz]} Baez, J.: The n-category cafe: Quivering with Excitement. Blog  May 4, 2012.
   http://golem.ph.utexas.edu/category/2012/05/quivering$\_$with$\_$excitement.html
\item{[Be]} Beilinson, A.A.: Coherent sheaves on $\Bbb P^n$ and problems in linear algebra.
    Funktsional. Anal. i Prilozhen 12.3 (1978), 68-29.
\item{[CR]} Caldero, Ph., Reineke, M.: On the quiver Grassmannian in the acyclic csae.
  J. Pure Appl. Algebra 212 (2008), 2369-2380.
    arXiv:math/0611074
\item{[D1]} Drozd, Yu.A.: Tame and wild matrix problems,
  in: Representation Theory II, Lecture Notes in Math., Vol. 832
  (1980), 242-258
\item{[D2]} Drozd, Yu.A.: Wild Algebras are Controlled Wild. Presentation ICRA XII, (2007). 
  http://icra12.mat.uni.torun.pl/lectures/Drozd-Torun.pdf 
\item{[Ha]}Harris, J.:  Lectures on rationally connected varieties. \newline
   http://mat.uab.es/\~kock/RLN/rcv.pdf.
\item{[L]} Le Bruyn, L.: Quiver Grassmannians can be anything. Blog, May 2, 2012.\newline 
    http://www.neverendingbooks.org/quiver-grassmannians-can-be-anything.
\item{[N]} Nakajima, H.; Quiver varieties and
   cluster algebras, http://arxiv.org/abs/0905.0002v5. 
\item{[Q]} Qin, F.: Quantum Cluster Variables via Serre Polynomials. \newline
    http://arxiv.org/abs/1004.4171.
\item{[Re]} Reineke, M.: Every projective variety is a quiver Grassmannian.
   Algebras and Representation Theory (to appear). arXiv:1204.5730.
\item{[Ri]}Ringel, C. M.: The Auslander bijections: How morphisms are determined by
  modules.  arXiv:1301.1251.
\item{[V]} Vakil, R.: 
   Murphy's Law in algebraic geometry: Badly-behaved deformation spaces.
   Invent. math. 164 (2006), 569-590.

	\bigskip\bigskip
{\rmk 
e-mail: \ttk ringel\@math.uni-bielefeld.de \par}

\bye